\documentclass[10pt,twoside]{article}
\usepackage{graphicx}
\usepackage{amsmath}
\usepackage{amssymb}
\usepackage{Latex-document}
\def\volumeno{I}
\title{\bf  Galois Representations\thanks{The work on this article was
partially supported by NSF Grant DMS-9702885.} \vskip 6mm}
\author{R. Taylor\vspace*{-0.5cm}\thanks{Department of Mathematics, Harvard
University, 1 Oxford St., Cambridge, MA 02138, USA. E-mail:
rtaylor@math.harvard.edu}}
\date{\vspace{-8mm}}

\newtheorem{thm}{Theorem}[section]
\newtheorem{prop}[thm]{Proposition}

\newtheorem{cor}[thm]{Corollary}
\newtheorem{conj}[thm]{Conjecture}

\newtheorem{conjecture}[thm]{Conjecture}

\newcommand{\onto}{\rightarrow \hspace{-.86em} \rightarrow}

\newcommand{\ra}{\rightarrow}
\newcommand{\lra}{\longrightarrow}
\newcommand{\la}{\leftarrow}

\newcommand{\into}{\hookrightarrow}

\newcommand{\iso}{\stackrel{\sim}{\ra}}
\newcommand{\liso}{\stackrel{\sim}{\lra}}

\newlength{\ownl}

\newcommand{\norm}{{\mbox{\bf N}}}
\newcommand{\ndiv}{{\mbox{$\not| $}}}

\newcommand{\Art}{{\operatorname{Art}\,}}

\newcommand{\End}{{\operatorname{End}\,}}

\newcommand{\Frob}{{\operatorname{Frob}}}
\newcommand{\Gal}{{\operatorname{Gal}\,}}

\newcommand{\Hom}{{\operatorname{Hom}\,}}

\newcommand{\Ind}{{\operatorname{Ind}\,}}

\newcommand{\Ree}{{\operatorname{Re}\,}}

\newcommand{\WD}{{\operatorname{WD}}}
\newcommand{\HT}{{\operatorname{HT}}}

\newcommand{\ad}{{\operatorname{ad}\,}}

\newcommand{\rec}{{\operatorname{rec}}}

\newcommand{\sol}{{\operatorname{sol}}}
\newcommand{\ab}{{\operatorname{ab}}}

\newcommand{\et}{{\operatorname{et}}}

\newcommand{\semis}{{\operatorname{ss}}}

\newcommand{\A}{{\mathbb{A}}}

\newcommand{\C}{{\mathbb{C}}}

\newcommand{\E}{{\mathbb{E}}}
\newcommand{\F}{{\mathbb{F}}}
\newcommand{\G}{{\mathbb{G}}}

\newcommand{\Q}{{\mathbb{Q}}}
\newcommand{\R}{{\mathbb{R}}}

\newcommand{\Z}{{\mathbb{Z}}}

\newcommand{\CA}{{\cal{A}}}

\newcommand{\CO}{{\cal{O}}}

\newcommand{\CR}{{\cal{R}}}

\newcommand{\gl}{{\mathfrak{gl}}}
\newcommand{\gm}{{\mathfrak{m}}}

\newcommand{\gz}{{\mathfrak{z}}}

\newcommand{\Fbar}{\overline{{\F}}}

\newcommand{\Qbar}{\overline{{\Q}}}
\newcommand{\barR}{\overline{{R}}}

\newcommand{\hatZ}{{\widehat{\Z}}}

\renewcommand{\thesection}{\arabic{section}}
\renewcommand{\thesubsection}{\thesection.\arabic{subsection}}
\def\Section#1{\section{\hskip -1em . \hskip 0.6em #1}}

\begin{document}
\maketitle

\thispagestyle{first} \setcounter{page}{449}

\begin{abstract}\vskip 3mm
In the first part of this paper we try to explain to a general mathematical
audience some of the remarkable web of conjectures linking representations of
Galois groups with algebraic geometry, complex analysis and discrete subgroups
of Lie groups. In the second part we briefly review some limited recent
progress on these conjectures.
\vskip 4.5mm

\noindent {\bf 2000 Mathematics Subject Classification:} 11F80.

\noindent {\bf Keywords and Phrases:} Galois representations,
$L$-function, Automorphic forms.
\end{abstract}

\vskip 12mm

\section*{Introduction}

\vskip-5mm \hspace{5mm}

The organisers requested a talk which would both be a colloquium style talk
understandable to a wide spectrum of mathematicians and one which would survey
the recent developments in the subject. I have found it hard to meet both
desiderata, and have opted to concentrate on the former. Thus the first three
sections of this paper contain a simple presentation of a web of deep
conjectures connecting Galois representations to algebraic geometry, complex
analysis and discrete subgroups of Lie groups. This will be of no interest
to the specialist. My hope is that the result is not too banal and that it
will give the non-specialist some idea of what motivates work in this area.
I should stress that nothing I write here is original.
In the final section I briefly review
some of what is known about these conjectures and {\em very briefly} mention
some of the available techniques. I also mention two questions which lie
outside the topic we are discussing, but which would have important implications
for it. Maybe someone can make progress on them?

Due to lack of space much of this article is too abbreviated. A somewhat
expanded version is available on my website {\tt www.math.harvard.edu/\~{}rtaylor}
and will hopefully be published elsewhere.

\Section{Galois representations}\label{s1}

\vskip-5mm \hspace{5mm}

We will let $\Q$ denote the field of rational numbers and $\Qbar$
denote the field of algebraic numbers, the algebraic closure of
$\Q$. We will also let $G_\Q$ denote the group of automorphisms of
$\Qbar$, that is $\Gal(\Qbar/\Q)$, the absolute Galois group of
$\Q$. Although it is not the simplest it is arguably the most natural
Galois group to study. An important technical point is that $G_\Q$ is
naturally a (profinite) topological group, a basis of open neighbourhoods of the
identity being given by the subgroups $\Gal(\Qbar/K)$ as $K$ runs over
subextensions of $\Qbar/\Q$ which are finite over $\Q$.

To my mind the Galois theory of $\Q$ is most interesting when one looks not
only at $G_\Q$ as an abstract (topological) group, but as a group with certain
additional structures associated to the prime numbers. I will now
briefly describe these structures.

For each prime number $p$ we may define an absolute value $| \,\,\,
|_p$ on $\Q$ by setting
\[ | \alpha |_p = p^{-r} \]
if $\alpha = p^r a/b$ with $a$ and $b$ integers coprime to $p$. If we
complete $\Q$ with respect to this absolute value we obtain the field
of $p$-adic numbers $\Q_p$, a totally disconnected, locally compact
topological field. We will write $G_{\Q_p}$ for its absolute Galois
group, $\Gal(\Qbar_p/\Q_p)$. The absolute value $|\,\,\,|_p$ has a
unique extension to an absolute value on $\Qbar_p$ and $G_{\Q_p}$ is
identified with the group of automorphisms of $\Qbar_p$ which preserve
$|\,\,\,|_p$, or equivalently the group of continuous automorphisms of
$\Qbar_p$. For each embedding $\Qbar \into \Qbar_p$ we obtain a closed
embedding $G_{\Q_p} \into G_\Q$ and as the embedding $\Qbar \into
\Qbar_p$ varies we obtain a conjugacy class of closed embeddings
$G_{\Q_p} \into G_\Q$. Slightly abusively, we shall consider $G_{\Q_p}$
a closed subgroup of $G_\Q$, suppressing the fact that the embedding is
only determined up to conjugacy.

This can be compared with the situation `at infinity'. Let
$|\,\,\,|_\infty$ denote the usual Archimedean absolute value on
$\Q$. The completion of $\Q$ with respect to $|\,\,\,|_\infty$ is the
field of real numbers $\R$ and its algebraic closure is $\C$ the field
of complex numbers. Each embedding $\Qbar \into \C$ gives rise to a
closed embedding
\[ \{ 1, c\} = G_\R = \Gal(\C/\R) \into G_\Q. \]
As the embedding $\Qbar \into \C$ varies one obtains a conjugacy class
of elements $c \in G_\Q$ of order $2$, which we refer to as complex
conjugations.

There are however many important differences between the case of
finite places (i.e. primes) and the infinite place $| \,\,\,
|_\infty$. For instance $\Qbar_p/\Q_p$ is an infinite extension and
$\Qbar_p$ is not complete. We will denote its completion by
$\C_p$. The Galois group $G_{\Q_p}$ acts on $\C_p$ and is in fact the
group of continuous automorphisms of $\C_p$.

The elements of $\Q_p$ (resp. $\Qbar_p$) with absolute
value less than or equal to $1$, form a closed subring $\Z_p$
(resp. $\CO_{\Qbar_p}$). These rings are local
with maximal ideals $p\Z_p$ (resp. $\gm_{\Qbar_p}$)
consisting of the elements with absolute value
strictly less than $1$. The field $\CO_{\Qbar_p}/
\gm_{\Qbar_p}$ is an algebraic closure of the
finite field with $p$ elements $\F_p = \Z_p/p\Z_p$, and we will
denote it by $\Fbar_p$. Thus we obtain a continuous map
\[ G_{\Q_p} \lra G_{\F_p} \]
which is surjective. Its kernel is called the inertia subgroup of
$G_{\Q_p}$ and is denoted by $I_{\Q_p}$. The group $G_{\F_p}$ is
procyclic and has a canonical generator called the (geometric)
Frobenius element and defined by
\[ \Frob_p^{-1} (x) = x^p. \]
In many circumstances it is technically convenient to replace $G_{\Q_p}$ by
a dense subgroup $W_{\Q_p}$, which is referred to as the Weil group of
$\Q_p$ and which is defined as the subgroup of $\sigma \in G_{\Q_p}$
such that $\sigma$ maps to
\[ \Frob_p^\Z \subset G_{\F_p}. \]
We endow $W_{\Q_p}$ with a topology by decreeing that $I_{\Q_p}$ with
its usual topology should be an open subgroup of $W_{\Q_p}$.

We will take a moment to describe some of the finer structure of
$I_{\Q_p}$ which we will need for technical purposes later. First of
all there is a (not quite canonical) continuous surjection
\[ I_{\Q_p} \onto \prod_{l \neq p} \Z_l \]
such that
\[ t(\Frob_p \sigma \Frob_p^{-1}) = p^{-1}t(\sigma) \]
for all $\sigma \in I_{\Q_p}$. The kernel of $t$ is a pro-$p$-group
called the wild inertia group. The fixed field $\Qbar_p^{\ker t}$ is
obtained by adjoining $\sqrt[n]{p}$ to $\Qbar_p^{I_{\Q_p}}$ for all
$n$ coprime to $p$ and
\[ \sigma \sqrt[n]{p} = \zeta_n^{t(\sigma)} \sqrt[n]{p}, \]
for some primitive $n^{th}$-root of unity $\zeta_n$ (independent of
$\sigma$, but dependent on $t$).

In my opinion the most interesting question about $G_\Q$ is to
describe it together with the distinguished subgroups $G_\R$,
$G_{\Q_p}$, $I_{\Q_p}$ and the distinguished elements $\Frob_p \in G_{\Q_p}/
I_{\Q_p}$.

I want to focus here on attempts to describe $G_\Q$  via its
representations. Perhaps the most obvious representations to consider
are those representations
\[ G_\Q \lra GL_n(\C) \]
with open kernel, and these so called Artin representations are
already very interesting. However one obtains a richer theory if one
considers representations
\[ G_\Q \lra GL_n(\Qbar_l) \]
which are continuous with respect to the $l$-adic topology on
$GL_n(\Qbar_l)$. We refer to these as {\em $l$-adic representations}.

One justification for considering $l$-adic representations is that
they arise naturally from geometry. Here are some examples of
$l$-adic representations.
\begin{enumerate}
\item A choice of embeddings $\Qbar \into \C$ and $\Qbar \into
\Qbar_l$ establishes a bijection between isomorphism classes of Artin
representations and isomorphism classes of $l$-adic representations
with open kernel. Thus Artin representations are a special case of
$l$-adic representations: those with finite image.

\item There is a unique character
\[ \chi_l:G_\Q \lra \Z_l^\times \subset \Qbar_l^\times \]
such that
\[ \sigma \zeta = \zeta^{\chi_l(\sigma)} \]
for all $l$-power roots of unity $\zeta$. This is called the $l$-adic
{\em cyclotomic character}.

\item If $X/\Q$ is a smooth projective variety (and we choose an embedding
$\Qbar \subset \C$) then the natural action of $G_\Q$ on the cohomology
\[ H^i(X(\C),\Qbar_l) \cong H^i_\et(X \times_\Q \Qbar, \Qbar_l) \]
is an $l$-adic representation. For instance if $E/\Q$ is an elliptic
curve then we have the concrete description
\[ H^1_\et(E \times_\Q \Qbar, \Qbar_l) \cong \Hom_{\Z_l}(\lim_{\la r}
E[l^r](\Qbar) , \Qbar_l) \cong \Qbar_l^2, \]
where $E[l^r]$ denotes the $l^r$-torsion points on $E$.
 We will write $H^i(X(\C),\Qbar_l(j))$
for the twist
\[ H^i(X(\C),\Qbar_l) \otimes \chi_l^{j}. \]
\end{enumerate}

Before discussing $l$-adic representations of $G_\Q$ further, let us
take a moment to look at $l$-adic representations of $G_{\Q_p}$. The
cases $l\neq p$ and $l=p$ are very different. Consider first the much
easier case $l \neq p$. Here $l$-adic representations of $G_{\Q_p}$
are not much different from representations of $W_{\Q_p}$ with open
kernel. More precisely define a {\em WD-representation} of $W_{\Q_p}$ over a
field $E$ to be a pair
\[ r:W_{\Q_p} \lra GL(V) \]
and
\[ N \in \End(V), \]
where $V$ is a finite dimensional $E$-vector space, $r$ is a
representation with open kernel and $N$ is a nilpotent endomorphism
which satisfies
\[ r(\phi)Nr(\phi^{-1}) = p^{-1}N \]
for every lift $\phi \in W_{\Q_p}$ of $\Frob_p$. The key point here is that
there is no reference to a topology on $E$, indeed no assumption that $E$ is
a topological field. Given $r$ there are
up to isomorphism only finitely many choices for the pair $(r,N)$ and
these can be explicitly listed without difficulty. A WD-representation
$(r,N)$ is called {\em unramified} if $N=0$ and $r(I_{\Q_p}) = \{ 1 \}$. It
is called Frobenius semi-simple if $r$ is semi-simple. Any
WD-representation $(r,N)$ has a canonical Frobenius
semi-simplification $(r,N)^\semis$ (see \cite{tate:ntb}).
In the case that $E=\Qbar_l$ we call
$(r,N)$ $l$-integral if all the eigenvalues of $r(\phi)$ have
absolute value $1$. This is independent of the choice of Frobenius
lift $\phi$.

If $l \neq p$, then there is an equivalence of categories between $l$-integral
WD-representations of $W_{\Q_p}$ over $\Qbar_l$ and $l$-adic
representations of $G_{\Q_p}$. To describe it choose a Frobenius lift
$\phi \in W_{\Q_p}$ and a surjection $t_l:I_{\Q_p} \onto \Z_l$. Up to
natural isomorphism the equivalence does not depend on these
choices. We associate to an $l$-integral WD-representation $(r,N)$ the unique
$l$-adic representation sending
\[ \phi^n \sigma \longmapsto r(\phi^n \sigma) \exp(t_l(\sigma)N) \]
for all $n \in \Z$ and $\sigma \in I_{\Q_p}$. The key point is
Grothendieck's observation that for $l \neq p$ any $l$-adic
representation of $G_{\Q_p}$ must be
trivial on some open subgroup of the wild inertia group. We will write
$\WD_p(R)$ for the WD-representation associated to an $l$-adic
representation $R$. Note that $\WD_p(R)$ is unramified if and only if
$R(I_p)=\{ 1\}$. In this case we call $R$ {\em unramified}.

The case $l=p$ is much more complicated because there are many more
$p$-adic representations of $G_{\Q_p}$. These have been extensively
studied by Fontaine and his co-workers. They
single out certain $p$-adic representations which they call {\em de Rham}
representations. I will not recall the somewhat involved definition
here (see however \cite{derham} and \cite{semisimp}), but note that `most' $p$-adic representations
of $G_{\Q_p}$ are not de Rham. To any de Rham representation $R$ of
$G_{\Q_p}$ on a $\Qbar_p$-vector space $V$ they associate the following.
\begin{enumerate}
\item A WD-representation $\WD_p(R)$ of $W_{\Q_p}$ over $\Qbar_p$ (see
\cite{berger} and \cite{fpr}).
\item A multiset $\HT(R)$ of $\dim V$ integers, called the Hodge-Tate
  numbers of $R$. The multiplicity of $i$ in $\HT(R)$ is
\[ \dim_{\Qbar_p} ( V \otimes_{\Q_p} \C_p(i))^{G_{\Q_p}}, \]
where $\C_p(i)$ denotes $\C_p$ with $G_{\Q_p}$-action $\chi_p(\sigma)^i$
times the usual (Galois) action on $\C_p$.
\end{enumerate}

We now return to the global situation (i.e. to the study of
$G_{\Q}$). The $l$-adic representations of $G_\Q$ that arise `in
nature', by
which I mean `from geometry', have a number of very special properties
which I will now list. Let $R:G_\Q \lra GL(V)$ be a subquotient of
$H^i(X(\C),\Qbar_l(j))$ for some smooth projective variety $X/\Q$ and
some
integers $i \geq 0$ and $j$.

\begin{enumerate}
\item (Grothendieck) The representation $R$ is unramified at all but finitely
many primes $p$.

\item (Fontaine, Messing, Faltings, Kato, Tsuji, de Jong, see e.g.
\cite{illu}, \cite{bert}) The representation
$R$ is de Rham in the sense that its restriction to $G_{\Q_l}$ is de Rham.

\item (Deligne, \cite{deligne:w1}) The representation $R$ is {\em pure} of weight $w=i-2j$ in
  the following sense. There is a finite set of primes $S$, such that
  for $p \not\in S$, the representation $R$ is unramified at $p$ and
for every eigenvalue $\alpha$ of $R(\Frob_p)$ and every
  embedding $\iota: \Qbar_l \into \C$
\[ | \iota\alpha |_\infty^2 = p^w. \]
In particular $\alpha$ is algebraic (i.e. $\alpha \in \Qbar$).
\end{enumerate}

An amazing conjecture of Fontaine and Mazur (see
\cite{fontaine:arbeit} and \cite{fm}) asserts that any
irreducible $l$-adic representation of $G_\Q$ satisfying the first two
of these properties arises from geometry in the above sense and so in
particular also satisfies the third property.

\begin{conjecture}[{\rm Fontaine-Mazur}]\label{cfm} Suppose that
\[ R:G_\Q \lra GL(V) \]
is an irreducible $l$-adic representation which is unramified at all
but finitely many primes and with $R|_{G_{\Q_l}}$ de Rham. Then there
is a smooth projective variety $X/\Q$ and integers $i \geq 0$ and $j$
such that $V$ is a subquotient of $H^i(X(\C),\Qbar_l(j))$. In
particular $R$ is pure of some weight $w\in \Z$.
\end{conjecture}

We will discuss the evidence for this conjecture later. We will call
an $l$-adic representation satisfying the conclusion of this
conjecture {\em geometric}.

Algebraic geometers have formulated some very precise conjectures about
the action of $G_\Q$ on the cohomology of varieties. We don't have the space
here to discuss these in general, but we will formulate, in an as algebraic
a way as possible, some of their conjectures.

\begin{conjecture}[{\rm Tate}]\label{cgeo} Suppose that $X/\Q$ is a smooth projective
variety. Then there is a decomposition
\[ H^i(X(\C), \Qbar) = \bigoplus_j M_j \]
with the following properties.
\begin{enumerate}
\item For each prime $l$ and for each embedding $\iota:\Qbar \into \Qbar_l,$
$M_j \otimes_{\Qbar, \iota} \Qbar_l$ is an irreducible subrepresentation of
$H^i(X(\C),\Qbar_l)$.
\item For all indices $j$ and for all primes $p$ there is a WD-representation
$\WD_p(M_j)$ of $W_{\Q_p}$ over $\Qbar$ such that
\[ \WD_p(M_j) \otimes_{\Qbar,\iota} \Qbar_l \cong \WD_p(M_j
\otimes_{\Qbar,\iota} \Qbar_l) \]
for all primes $l$ and all embeddings $\iota:\Qbar \into \Qbar_l$.
\item There is a multiset of integers $\HT(M_j)$ such that
\begin{enumerate}
\item for all primes $l$ and all embeddings $\iota:\Qbar \into \Qbar_l$
\[ \HT(M_j \otimes_{\Qbar,\iota} \Qbar_l) = HT(M_j)\]
\item and for all $\iota:\Qbar \into \C$
\[ \dim_\C ((M_j \otimes_{\Qbar,\iota} \C) \cap H^{a,i-a}(X(\C),\C)) \]
is the multiplicity of $a$ in $HT(M_j)$.
\end{enumerate} \end{enumerate} \end{conjecture}

If one considers the whole of $ H^i(X(\C), \Qbar)$ rather than its pieces
$M_j$, then part 2. is known to hold up to Frobenius semisimplification for
all but finitely many $p$ and part 3. is known to hold (see \cite{illu}). It follows from a
theorem of Faltings \cite{falt:mor} that the whole conjecture is true for $H^1$ of an
abelian variety. The putative constituents $M_j$ are one incarnation of what
people call `pure motives'.

If one believes conjectures \ref{cfm} and \ref{cgeo} then `geometric' $l$-adic
representations should come in compatible families as $l$ varies. There are
many ways to make precise the notion of such a compatible family. Here is one.

By a {\em weakly compatible system of $l$-adic representations} $\CR=\{
R_{l,\iota} \}$ we shall mean a collection of semi-simple $l$-adic
representations
\[ R_{l,\iota}: G_{\Q} \lra GL(V \otimes_{\Qbar,\iota} \Qbar_l), \]
one for each pair $(l,\iota)$, where $l$ is a prime and $\iota:\Qbar
\into \Qbar_{l}$, which satisfy the following conditions.
\begin{itemize}
\item There is a multiset of integers $\HT(\CR)$ such that for each prime $l$
and each embedding $\iota:\Qbar \into \Qbar_l$ the restriction
$R_{l,\iota}|_{G_{\Q_l}}$ is de Rham and
$\HT(R_{l,\iota}|_{G_{\Q_l}}) = \HT(\CR)$.
\item There is a finite set of primes $S$ such that if $p \not\in S$
then $\WD_p(R_{l,\iota})$ is unramified for all $l$ and $\iota$.
\item For all but finitely many primes $p$ there is a Frobenius
semi-simple WD-representation $\WD_p(\CR)$ over $\Qbar$ such that for all primes $l \neq p$
and for all $\iota$ we have
\[ \WD_p(R_{l,\iota})^\semis \sim \WD_p(\CR). \]
\end{itemize}
We make the following subsidiary definitions.
\begin{itemize}
    \item We call $\CR$ {\em strongly compatible} if the last condition
    (the existence of $\WD_p(\CR)$) holds for all primes $p$.
    \item We call $\CR$ {\em irreducible} if each $R_{l,\iota}$ is
    irreducible.
    \item We call $\CR$ {\em pure} of weight $w \in \Z$, if for all but
    finitely many $p$ and for all eigenvalues $\alpha$ of
    $r_p(\Frob_p)$, where $\WD_p(\CR)=(r_p,N_p)$, we have
    \[ |\iota \alpha|_\infty^2 = p^w \]
    for all embeddings $\iota:\Qbar \into \C$.
    \item We call $\CR$ {\em geometric} if there is a smooth projective variety
    $X/\Q$ and integers $i \geq 0$ and $j$ and a subspace
    \[ W \subset H^i(X(\C),\Qbar) \]
    such that for all $l$ and $\iota$, $W \otimes_{\Qbar, \iota} \Qbar_l$ is
    $G_\Q$ invariant and realises $R_{l,\iota}$.
\end{itemize}

Conjectures \ref{cfm} and \ref{cgeo} lead one to make the following
conjecture.

\begin{conjecture}\label{ccs}\begin{enumerate}
\item If $R:G_\Q \ra GL_n(\Qbar_l)$ is a continuous semi-simple de Rham
representation unramified at all but finitely many primes then $R$ is part of
a weakly compatible system.
\item Any weakly compatible system is strongly compatible.
\item Any irreducible weakly compatible system $\CR$ is geometric and pure of weight
$(2/\dim \CR) \sum_{h \in \HT(\CR)} h$.
\end{enumerate} \end{conjecture}

A famous theorem of Cebotarev asserts that if $K/\Q$ is any Galois extension
in which all but finitely many primes are unramified (i.e. for all but finitely
many primes $p$ the image of $I_{\Q_p}$ in $\Gal(K/\Q)$ is trivial) then the
Frobenius elements at unramified primes $\Frob_p \in \Gal(K/\Q)$ are dense in
$\Gal(K/\Q)$. It follows that an irreducible weakly compatible system $\CR$ is
uniquely determined by $\WD_p(\CR)$ for all but finitely many $p$ and hence by one $R_{l,\iota}$.

Conjectures \ref{cfm} and \ref{ccs} are known for one dimensional representations,
in which case they have purely algebraic proofs based on class
field theory (see \cite{serre:alr}). Otherwise
only fragmentary cases have been proved, where amazingly the arguments
are extremely indirect involving sophisticated analysis and geometry. We
will come back to this later.

\Section{{\boldmath $L$}-functions}\label{s2}

\vskip-5mm \hspace{5mm}

$L$-functions are certain Dirichlet series
\[ \sum_{n=1}^\infty a_n/n^s \]
which play an important role in number theory. A full discussion
of the role of $L$-functions in number theory is beyond the scope
of this talk. The simplest example of an $L$-function is the
Riemann zeta function
\[ \zeta(s)=\sum_{n=1}^\infty 1/n^s. \]
It converges to a holomorphic function in the half plane $\Ree s > 1$ and in
this region of convergence it can also be expressed as a convergent infinite
product over the prime numbers
\[ \zeta(s) = \prod_p (1-1/p^s)^{-1}. \]
This is called an {\em Euler product} and the individual factors are called Euler
factors. Lying deeper is the fact that $\zeta(s)$ has meromorphic continuation
to the whole complex plane, with only one pole: a simple pole at $s=1$.
Moreover if we set
\[ Z(s)=\pi^{-s/2} \Gamma(s/2) \zeta(s) \]
then $Z$ satisfies the functional equation
\[ Z(1-s)=Z(s). \]
Encoded in the Riemann zeta function is lots of deep arithmetic information.
For instance the location of the zeros of $\zeta(s)$ is intimately connected
with the distribution of prime numbers. Moreover its special values at negative
integers (where it is only defined by analytic continuation) turn out to be
rational numbers encoding deep arithmetic information about the cyclotomic
fields $\Q(e^{2 \pi \sqrt{-1}/p})$.

Another celebrated example is the $L$-function of an elliptic
curve $E$:
\[ y^2=x^3+ax+b. \]
In this case the $L$-function is defined as an Euler product
(converging in $\Ree s>3/2$)
\[ L(E,s) = \prod_p L_p(E,p^{-s}), \]
where $L_p(E,X)$ is a rational function, and for all but finitely many $p$
\[ L_p(E,X) = (1-a_p(E)X+pX^2)^{-1}, \]
with $p-a_p(E)$ being the number of solutions to the congruence
\[ y^2 \equiv x^3+ax+b \bmod p \]
in $\F_p^2$. It has recently been proved \cite{bcdt} that $L(E,s)$ can be
continued to an entire function, which satisfies a functional equation
\[ (2 \pi)^{-s} \Gamma(s) L(E,s)= \pm N(E)^{1-s} (2 \pi)^{s-2}
\Gamma(2-s) L(E,2-s), \] for some explicit positive integer
$N(E)$. A remarkable conjecture of Birch and Swinnerton-Dyer
\cite{bsd} predicts that $y^2=x^3+ax+b$ has infinitely many
rational solutions if and only if $L(E,1)=0$. Again we point out
that it is the behaviour of the $L$-function at a point where it
is only defined by analytic continuation, which is governing the
arithmetic of $E$. This conjecture has been proved (see
\cite{koly}) when $L(E,s)$ has at most a simple zero at $s=1$.

One general setting in which one can define $L$-functions is
$l$-adic representations. Let us look first at the local setting.
If $(r,N)$ is a WD-representation of $W_{\Q_p}$ on an $E$-vector
space $V$, where $E$ is an algebraically closed field of
characteristic zero, we define a local L-factor
\[ L((r,N),X) = \det(1-X \Frob_p)|_{V^{I_{\Q_p},N=0}}^{-1} \in E(X). \]
($V^{I_{\Q_p},N=0}$ is the subspace of $V$ where $I_{\Q_p}$ acts trivially
and $N=0$.)
One can also associate to $(r,N)$ a conductor $f(r,N) \in \Z_{\geq 0}$,
which measures how deeply into $I_{\Q_p}$ the WD-representation $(r,N)$ is
nontrivial, and a local epsilon factor $\epsilon((r,N),\Psi_p) \in E$, which
also depends on the choice of a non-trivial character $\Psi_p:\Q_p
\ra E^\times$ with open kernel. (See \cite{tate:ntb}.)

If $R:G_\Q \ra GL(V)$ is an $l$-adic representation of $G_\Q$
which is de Rham at $l$ and pure of some weight $w \in \Z$, and if
$\iota:\Qbar_l \into \C$ we will define an $L$-function
\[ L(\iota R,s) = \prod_p L(\iota \WD_p(R),p^{-s}), \]
which will converge to a holomorphic function in $\Ree s> 1+w/2$. For example
\[ L(1,s) = \zeta(s) \]
and if $E/\Q$ is an elliptic curve then
\[ L(\iota H^1(E(\C),\Qbar_l), s) = L(E,s) \]
(for any $\iota$). Note the useful formulae
\[ L(\iota(R_1 \oplus R_2),s) = L(\iota R_1,s)L(\iota R_2,s)
\,\,\,\,\,\,\,\,\,\, {\rm and } \,\,\,\,\,\,\,\,\,\,
L(\iota (R \otimes \chi_l^r),s) = L(\iota R, s+r). \]
Also note that $L(\iota R,s)$ determines $L(\WD_p(R),X)$ for all $p$ and
hence $\WD_p(R)^\semis$ for all but finitely many $p$. Hence by the Cebotarev
density theorem $L(\iota R,s)$ determines $R$ (up to semisimplification).

Write $m^R_i$ for the multiplicity of an integer $i$ in $\HT(R)$ and, if
$w/2 \in \Z$, define $m_{w/2,\pm}^R \in (1/2)\Z$ by:
\[ \begin{array}{rcl} m_{w/2,+}^R+m_{w/2,-}^{R} &=& m^R_{w/2} \\ \\
m_{w/2,+}^R-m_{w/2,-}^{R} &=&
(-1)^{w/2}(\dim V^{c=1} - \dim V^{c=-1}). \end{array} \]
Assume that $m^R_{w/2,\pm}$ are integers, i.e. that $m^R_{w/2} \equiv \dim V
\bmod 2$. Then we can define a $\Gamma$-factor, $\Gamma(R,s)$, which is a
product of functions $\pi^{-(s+a)/2}\Gamma((s+a)/2)$ as $a$ runs over a set
of integers depending only on the numbers $m^R_i$ and $m^R_{w/2,\pm}$. We can
also define an epsilon factor $\epsilon_\infty(R,\Psi_\infty) \in \C^\times$
which again only depends on $m^R_i$, $m^R_{w/2,\pm}$ and a non-trivial
character $\Psi_\infty:\R \ra \C^\times$. Set
\[ \Lambda(\iota R,s)= \Gamma(R,s) L(\iota R, s) \]
and
\[ N(R) = \prod_p p^{f(\WD_p(R))} \]
(which makes sense as $f(\WD_p(R))=0$ for all but finitely many $p$) and
\[ \epsilon(\iota R) = \epsilon_{\infty}(R,e^{2 \pi \sqrt{-1} x})
\prod_p \iota \epsilon(\WD_p(R), \Psi_p), \]
where $\iota \Psi_p(x) = e^{-2 \pi \sqrt{-1} x}$.

The following conjecture is a combination of conjecture \ref{cfm} and
conjectures which have become standard.
\begin{conjecture}\label{cfe} Suppose that $R$ is an irreducible $l$-adic
representation of $G_\Q$ which is de Rham and pure of weight $w \in \Z$.
Then $m_p^R=m_{w-p}^R$ for all $p$, so that $m_{w/2} \equiv
\dim V \bmod 2$. Moreover the following should hold.
\begin{enumerate}
\item $L(\iota R,s)$ extends to an entire function, except for a single simple
pole if $R=\chi_l^{-w/2}$.
\item $\Lambda(\iota R,s)$ is bounded in vertical strips $\sigma_0
\leq \Ree s \leq \sigma_1$.
\item $\Lambda(\iota R,s)= \epsilon(\iota R) N(R)^{-s}
\Lambda(\iota R^\vee,1-s)$.
\end{enumerate} \end{conjecture}

It is tempting to believe that something like properties 1., 2. and 3. should
characterise those Euler products which arise from $l$-adic representations.
We will discuss a more precise conjecture along these lines in the
next section. Why Galois representations should be {\bf the} source of Euler
products with good functional equations seems a complete mystery.

\Section{Automorphic forms}\label{s3}

\vskip-5mm \hspace{5mm}

Automorphic forms may be thought of as
certain smooth functions on the quotient $GL_n(\Z)\backslash GL_n(\R)$. We need
several preliminaries before we can make a precise definition.

Let $\hatZ$ denote the profinite completion of $\Z$, i.e.
\[ \hatZ = \lim_{\la N} \Z / N\Z = \prod_p \Z_p, \]
a topological ring. Also let $\A^\infty$ denote the topological ring of
finite adeles
\[ \A^\infty = \hatZ \otimes_\Z \Q, \]
where $\hatZ$ is an open subring with its usual topology. As an abstract
ring, $\A^\infty$ is the subring of $\prod_p \Q_p$ consisting of elements
$(x_p)$ with $x_p \in \Z_p$ for all but finitely many $p$. However the
topology is not the subspace topology. We define the topological ring of
adeles to be the product
\[ \A = \A^\infty \times \R. \]
Note that $\Q$ embeds diagonally as a discrete subring of $\A$ with
compact quotient
\[ \Q \backslash \A = \hatZ \times \Z \backslash \R. \]

We will be interested in $GL_n(\A)$, the locally compact topological group of $n
\times n$ invertible matrices with coefficients in $\A$. We remark
that the topology on $GL_n(\A)$ is the subspace topology resulting
from the closed embedding
\[  \begin{array}{rcl} GL_n(\A) &\into & M_n(\A) \times M_n(\A) \\
                         g & \mapsto & (g, g^{-1}). \end{array} \]
$GL_n(\Q)$ is a discrete subgroup of $GL_n(\A)$ and the quotient
$GL_n(\Q) \backslash GL_n(\A)$ has finite volume. If $U \subset GL_n(\hatZ)$
is an open subgroup with $\det U = \hatZ^\times$, then
\[ GL_n(\Q) \backslash GL_n(\A) / U = (GL_n(\Q) \cap U) \backslash GL_n(\R). \]
Note that $GL_n(\Q) \cap U$ is a subgroup of $GL_n(\Z)$ of finite index.
Most of the statements we make
concerning $GL_n(\A)$ can be rephrased to involve only $GL_n(\R)$, but at the
expense of making them much more cumbersome. To achieve brevity (and because
it seems more natural) we have opted to use the language of adeles. We hope
that this extra abstraction will not be too confusing for the novice.

Before continuing our introduction of automorphic forms let us digress to
mention class field theory, which provides a concrete example of the
presentational advantages of the adelic language. It also implies
essentially all the conjectures we are considering in the case of one
dimensional Galois representations. Indeed this article is about the
search for a non-abelian analogue of class field theory. Class field theory
gives a concrete description of the abelianisation (maximal continuous
abelian quotient) $G_\Q^\ab$ of $G_\Q$ and $W_{\Q_p}^\ab$ of
$W_{\Q_p}$. First the local theory asserts that there is an
isomorphism
\[ \Art_p: \Q_p^\times \iso W_{\Q_p}^\ab \]
with various natural properties, including the facts that $\Art(\Z_p^\times)$
is the image of the inertia group $I_{\Q_p}$ in $W_{\Q_p}^\ab$, and that the
induced map
\[ \Q_p^\times / \Z_p^\times \lra W_{\Q_p}^\ab/I_{\Q_p} \subset G_{\F_p} \]
takes $p$ to the geometric Frobenius element $\Frob_p$. Secondly the global
theory asserts that there is an isomorphism
\[ \Art: \A^\times/\Q^\times \R^\times_{>0} \iso G_\Q^\ab \]
such that the restriction of $\Art$ to $\Q_p^\times$ coincides with the
composition of $\Art_p$ with the natural map $W_{\Q_p}^\ab \ra G_\Q^\ab$.
Thus $\Art$ is defined completely from a knowledge of the $\Art_p$ (and
the fact that $\Art$ takes $-1 \in \R^\times$ to complex conjugation) and
the reciprocity theorem of global class field theory can be thought of as a
determination of the kernel of $\prod_p \Art_p$.

We now return to our (extended) definition of automorphic forms. For
each partition $n=n_1+n_2$ let $N_{n_1,n_2}$ denote the subgroup of
$GL_n$ consisting of matrices of the form
\[ \left( \begin{array}{cc} I_{n_1} & * \\ 0 & I_{n_2} \end{array}
\right) . \]
Let $O(n) \subset GL_n(\R)$ denote the orthogonal subgroup.
Let $\gz_n$ denote the centre of the universal enveloping of $\gl_n$,
the complexified Lie algebra of $GL_n(\R)$ (i.e. $\gl_n=M_n(\C)$ with
$[X,Y]=XY-YX$).
Via the Harish-Chandra isomorphism (see for example \cite{hc}) we may
identify homomorphisms $\gz_n \ra \C$ with multisets of $n$ complex numbers.
We will write $\chi_H$ for the homomorphism corresponding to a
multiset $H$. Thus $\gz_n$ acts on the irreducible finite dimensional
$\gl_n$-module with highest weight $(a_1,...,a_n) \in \Z^n$ ($a_1 \geq ...
\geq a_n$) by $\chi_{\{ a_1+(n-1)/2,...,a_n+(1-n)/2 \} }$.

Fix such a multiset $H$ of cardinality $n$. The space of cusp
forms with infinitesimal character $H$, $\CA^\circ_H(GL_n(\Q) \backslash
GL_n(\A))$ is the space of smooth bounded functions
\[ f: GL_n(\Q) \backslash GL_n(\A) \lra \C \]
satisfying the following conditions.
\begin{enumerate}
\item ($K$-finiteness) The translates of $f$ under $GL_n(\hatZ) \times O(n)$
(where $O(n)$ denotes the orthogonal group) span a finite dimensional
vector space;
\item (Infinitesimal character $H$) If $z \in \gz_n$ then $zf=\chi_H(z) f$;
\item (Cuspidality) For each partition $n=n_1+n_2$,
\[ \int_{N_{n_1,n_2}(\Q) \backslash N_{n_1,n_2}(\A)} f(ug) du =0. \]
\end{enumerate}
Note that if $U \subset GL_n(\hatZ)$ is an open subgroup with $\det
U = \hatZ^\times$ then one may think of $\CA^\circ_H(GL_n(\Q) \backslash
GL_n(\A))^U$ as a space of functions on $(GL_n(\Q) \cap U)
\backslash GL_n(\R)$.

One would like to study $\CA^\circ_H(GL_n(\Q) \backslash GL_n(\A))$ as a
representation of $GL_n(\A)$, unfortunately it is not preserved by the
action of $GL_n(\R)$ (because the $K$-finiteness condition depends on the
choice of a maximal compact subgroup $O(n) \subset GL_n(\R)$). It does however
have an action of $GL_n(\A^\infty) \times O(n)$ and of $\gl_n$, which is
essentially as good. More precisely it is an admissible $GL_n(\A^\infty)
\times (\gl_n,O(n))$-module in the sense of \cite{flath}. In fact it is a
direct sum of irreducible, admissible $GL_n(\A^\infty) \times
(\gl_n,O(n))$-modules each occurring with multiplicity one. We will
(slightly abusively) refer to these irreducible constituents as
cuspidal automorphic representations of $GL_n(\A)$ with
infinitesimal character $H$.

$\CA^\circ_{\{ 0 \} }(\Q^\times \backslash \A^\times)$ is
just the space of locally constant functions on $\A^\times/\Q^\times
\R^\times_{>0}$ and so cuspidal automorphic representations of $GL_1(\A)$
with infinitesimal character $\{ 0 \}$, are just the (finite order) complex
valued characters of $\A^\times/\Q^\times \R^\times_{>0} \cong \hatZ^\times$,
i.e. Dirichlet characters. $\CA^\circ_{\{ s \} }(\Q^\times \backslash
\A^\times)$ is simply obtained from $\CA^\circ_{\{ 0 \} }(\Q^\times
\backslash \A^\times)$ by twisting by $|| \,\,\,||^s$, where $||\,\,\,||:
\A^\times/\Q^\times \ra \R^\times_{>0}$ is the
product of the absolute values $| \,\,\,|_x$. Thus in the case $n=1$
cuspidal automorphic representations are essentially Dirichlet characters.

The case $n=2$ is somewhat more representative. In this case we have
$\CA^\circ_{ \{ s,t \} }(GL_2(\Q) \backslash GL_2(\A))=(0)$ unless $s-t \in
i\R$, $s-t \in \Z$ or $s-t \in (-1,1)$. It is conjectured that the third
possibility can not arise unless $s=t$. Let us consider the case $s-t \in
\Z_{>0}$ a little further. If $s-t \in \Z_{>0}$ then it turns out that
the irreducible constituents of $\CA^\circ_{ \{ s,t \} }(GL_2(\Q)
\backslash GL_2(\A))$ are in bijection with the weight $1+s-t$ holomorphic
cusp forms on the upper half plane, which are normalised newforms (see
for example \cite{miy}). Thus in some sense
cuspidal automorphic representations are are also generalisations of
classical holomorphic normalised newforms.

Note that if $\psi$ is a character of $\A^\times/\Q^\times \R^\times_{>0}$
and if $\pi$ is a cuspidal automorphic representation of $GL_n(\A)$ with
infinitesimal character $H$ then $\pi \otimes (\psi \circ \det)$ is
also a cuspidal automorphic representation with infinitesimal
character $H$ and the contragredient (dual) $\pi^*$ of $\pi$ is a cuspidal automorphic
representation with infinitesimal character $-H=\{ -h:\,\, h \in H \}$.

One of the main questions in the theory of automorphic forms is to describe
the irreducible constituents of $\CA^\circ_H(GL_n(\Q) \backslash GL_n(\A))$.
If we are to do this we first need some description of all irreducible
admissible $GL_n(\A^\infty) \times (\gl_n,O(n))$-modules, and then we can try
to say which occur in $\CA^\circ_H(GL_n(\Q) \backslash GL_n(\A))$.

Just as a character $\psi: \A^\times \ra \C^\times$ can be factored as
\[ \psi = \psi_\infty \times \prod_p \psi_p \]
where $\psi_p:\Q_p^\times \ra \C^\times$ (resp. $\psi_\infty:\R^\times \ra
\C^\times$), so an irreducible, admissible $GL_n(\A^\infty) \times
(\gl_n,O(n))$-module can be factorised as a restricted tensor
product (see \cite{flath})
\[ \pi \cong {\bigotimes_x}' \pi_x, \]
where $\pi_{\infty}$ is an irreducible, admissible $(\gl_n,O(n))$-module
(see for example \cite{wall}), and each $\pi_p$ is an irreducible smooth
(i.e. stabilisers of vectors are open in $GL_n(\Q_p)$) representation of
$GL_n(\Q_p)$ with $\pi_p^{GL_n(\Z_p)} \neq (0)$ for all but finitely
many $p$. To the factors $\pi_x$ one can associate various invariants
(see \cite{jac}).
\begin{itemize}
    \item A central character $\psi_x:\Q_x^\times \ra \C^\times$.
    \item L-factors $L(\pi_p,X) \in \C(X)$.
    \item A $\Gamma$-factor $\Gamma(\pi_{\infty},s)$.
    \item Conductors $f(\pi_p) \in \Z_{\geq 0}$.
    \item For each non-trivial character $\Psi_x:\Q_x \ra
    \C^{\times}$ an epsilon factor $\epsilon(\pi_x,\Psi_x) \in
    \C^{\times}$.
\end{itemize}
(We also remark that $\gz_n$ acts via a character $\gz_n \ra \C$ on
any irreducible, admissible $(\gl_n,O(n))$-module $\pi_\infty$. This
character is called the infinitesimal character of $\pi_\infty$.)
Now we may attach to $\pi$
\begin{itemize}
\item a central character $\psi_\pi = \prod_x \psi_x: \A^\times \ra \C^\times$;
\item an $L$-function $L(\pi,s) = \prod_p L(\pi_p,p^{-s})$ (which may or may not
converge);
\item an extended $L$-function $\Lambda(\pi,s)=\Gamma(\pi_\infty,s)L(\pi,s)$;
\item a conductor $N(\pi) = \prod_p p^{f(\pi_p)}$ (which makes sense because
$f(\pi_p)=0$ when $\pi_p^{GL_n(\Z_p)} \neq (0)$);
\item and an epsilon constant $\epsilon(\pi)
= \prod_x \epsilon(\pi_x,\Psi_x) \in \C^{\times}$ where $\prod_x
\Psi_x: \A/\Q \ra \C^\times$ is any non-trivial character.
\end{itemize}

The following theorem and conjecture describe the (expected)
relationship between automorphic forms and $L$-functions with
Euler product and functional equation. We suppose $n>1$. A similar
theorem to theorem \ref{gjt} is true for $n=1$, except that
$L(\pi,s)$ may have one simple pole. In this case it was due to
Dirichlet. Conjecture \ref{cct} becomes vacuous if $n=1$.

\begin{thm}[{\rm Godement-Jacquet, \cite{gj}}]\label{gjt} Suppose that $\pi$ is an irreducible \linebreak
constituent of $\CA^\circ_H(GL_n(\Q) \backslash GL_n(\A))$ with $n>1$. Then $L(\pi,s)$
converges to a holomorphic function in some right half complex plane
$\Ree s > \sigma$ and can be continued to a holomorphic function on the whole
complex plane so that $\Lambda(\pi,s)$ is bounded in all vertical strips $\sigma_1 \geq \Ree s \geq
\sigma_2$. Moreover $\Lambda(\pi,s)$ satisfies the
functional equation
\[ \Lambda(\pi,s)=\epsilon(\pi) N(\pi)^{-s} \Lambda(\pi^*,1-s). \]
\end{thm}

\begin{conj}[{\rm Cogdell-Piatetski-Shapiro, \cite{cps1}}]\label{cct} Suppose that $\pi$ is an irreducible, admissible
$GL_n(\A^\infty) \times (\gl_n,O(n))$-module such that the central
character of $\pi$ is trivial on $\Q^\times$ and such that $L(\pi,s)$ converges
in some half plane. Suppose also that for all characters $\psi:
\A^\times/\Q^\times \R^\times_{>0} \ra \C^\times$ the $L$-function
$\Lambda (\pi
\otimes (\psi \circ \det),s)$ (which will then converge in some right
half plane) can be continued to a holomorphic function
on the entire complex plane, which is bounded in vertical strips and satisfies
a functional equation
\[ \Lambda(\pi \otimes (\psi \circ \det),s)=\epsilon(\pi \otimes
(\psi \circ \det)) N(\pi \otimes (\psi \circ \det))^{-s}
\Lambda(\pi^* \otimes (\psi^{-1} \circ \det),1-s). \]
($\Lambda (\pi^* \otimes (\psi^{-1} \circ \det),s)$ also automatically
converges in some right half plane.)
Then there is a partition $n=n_1+...+n_r$ and cuspidal automorphic
representations $\pi_i$ of $GL_{n_i}(\A)$ such that
\[ \Lambda(\pi,s) = \prod_{i=1}^r \Lambda(\pi_i,s). \]
\end{conj}

This conjecture is known to be true for $n=2$ (\cite{weil}, \cite{jl}) and
$n=3$ (\cite{jpss:gl3}). For $n>3$ a weaker form of
this conjecture involving twisting by higher dimensional automorphic
representations is known to hold (see \cite{cps1}, \cite{cps2}). These
results are called `converse theorems'.

The reason for us introducing automorphic forms is because of a putative
connection to Galois representations, which we will now discuss. But
first let us discuss the local situation. It has recently been
established (\cite{ht}, \cite{hen}, \cite{har:icm}) that there is a natural
bijection, $\rec_p$, from
irreducible smooth representations of $GL_n(\Q_p)$ to $n$-dimensional
Frobenius semi-simple WD-representations of $W_{\Q_p}$ over $\C$. The
key point here is that the bijection should be natural. We will not describe here
exactly what this means (instead we refer the reader to the
introduction to \cite{ht}). It does satisfy the following.
\begin{itemize}
\item $\psi_{\pi} \circ \Art_p^{-1} = \det \rec_p (\pi)$, where
$\psi_\pi$ is the central character of $\pi$.
\item $L(\rec_p(\pi),X)= L(\pi,X)$.
\item $f(\rec_p(\pi))=f(\pi)$.
\item $\epsilon(\rec_p(\pi), \Psi_p)=\epsilon(\pi,\Psi_p)$ for any
non-trivial character $\Psi_p:\Q_p \ra \C^\times$.
\end{itemize}
The existence of $\rec_p$ can be seen as a non-abelian
generalisation of local class field theory, as in the case $n=1$ we have
$\rec_p(\pi)=\pi \circ \Art_p^{-1}$.

Now suppose that $\iota:\Qbar_l \ra \C$ and that $R$ is a de Rham
$l$-adic representation of $G_\Q$ which is unramified at all but
finitely many primes. Using the local reciprocity map
$\rec_p$, we can associate to $R$ an irreducible, admissible $GL_n(\A^\infty)
\times (\gl_n,O(n))$-module
\[ \pi(\iota R)= \pi_{\infty}(R) \otimes \prod_p \rec_p^{-1}(\iota
\WD_p(R)), \]
where $\pi_{\infty}(R)$ is a tempered irreducible, admissible
$(\gl_n,O(n))$-module with infinitesimal character
$\HT(R|_{G_{\Q_l}})$ and with
$\Gamma(\pi_{\infty}(R),s)=\Gamma(R,s)$. The definition of
$\pi_\infty(R)$ depends only on the numbers $m^R_i$ and
$m^R_{w/2,\pm}$. Then we have the following conjectures.

\begin{conj}\label{cgfa} Suppose that $H$ is a multiset of $n$
{\em integers} and that $\pi$ is an
irreducible constituent of $\CA^\circ_H(GL_n(\Q) \backslash GL_n(\A))$.
Identify $\Qbar \subset \C$. Then each $\rec_p(\pi_p)$ can be defined over
$\Qbar$ and there is an irreducible geometric strongly compatible system of $l$-adic
representations $\CR$ such that $\HT(\CR)=H$ and
$\WD_p(\CR)^{\semis}=\rec_p (\pi_p)$ for all primes $p$. \end{conj}

\begin{conj}\label{cafg} Suppose that
\[ R:G_\Q \lra GL(V) \]
is an irreducible $l$-adic representation which is unramified at all but
finitely many primes and for which $R|_{G_{\Q_l}}$ is de Rham. Let $\iota:
\Qbar_l \ra \C$. Then $\pi(\iota R)$ is a cuspidal automorphic
representation of $GL_n(\A)$.
\end{conj}

These conjectures are essentially due to Langlands \cite{lang:conj}, except
we have used a precise formulation which follows Clozel \cite{claa} and we
have incorporated conjecture \ref{cfm} into conjecture \ref{cafg}.

Conjecture \ref{cafg} is probably the more mysterious of the two, as
only the case $n=1$ and fragmentary cases where $n=2$ are known. This
will be discussed further in the next section. Note the
similarity to the main theorem of global class field theory that
$\prod_p \Art_p:\A^\times \ra G_\Q^\ab$ has kernel $\Q^\times$.

The following theorem provides significant evidence for conjecture
\ref{cgfa}.

\begin{thm}[{\rm \cite{kotjams}, \cite{clihes}, \cite{ht}}]\label{kt} Suppose that
$H$ is multiset of $n$ {\em distinct} integers and that $\pi$ is an irreducible constituent of
$\CA^\circ_H(GL_n(\Q) \backslash GL_n(\A))$. Let $\iota:\Qbar_l \into \C$.
Suppose moreover that
$\pi^* \cong \pi \otimes (\psi \circ \det)$ for some character
$\psi:\A^\times/\Q^\times \ra \C^\times$, and that either $n \leq 2$ or for some prime $p$ the
representation $\pi_p$ is square integrable (i.e. $\rec_p(\pi_p)$ is indecomposable). Then there is a continuous
representation
\[ R_{l,\iota}: G_\Q \lra GL_n(\Qbar_l) \]
with the following properties.
\begin{enumerate}
\item $R_{l,\iota}$ is geometric and pure of weight $2/n \sum_{h \in H}
h$.
\item $R_{l,\iota}|_{G_{\Q_l}}$ is de Rham and $\HT(R_{l,\iota}|_{G_{\Q_l}})=
H$.
\item For any prime $p \neq l$ there is a representation $r_p:W_{\Q_p} \ra
GL_n(\Qbar_l)$ such that $\WD_p(R_{l,\iota})^\semis=(r_p,N_p)$ and $\rec_p(\pi_p)=
(\iota r_p, N_p')$.
\end{enumerate}
\end{thm}

This was established by finding the desired $l$-adic
representations in the cohomology of certain unitary group Shimura varieties.
It seems not unreasonable to hope that similar techniques might allow one to
improve many of the technical defects in the theorem.
However Clozel has stressed that in the cases where $H$ does not have distinct
elements or where $\pi^* \not\cong \pi \otimes(\psi \circ \det)$, there seems to be
no prospect of finding the desired $l$-adic representations in the cohomology
of Shimura varieties. It seems we need a new technique.

\Section{What do we know?} \label{s4}

\vskip-5mm \hspace{5mm}

Let us first summarise in a slightly less precise way the various conjectures
we have made, in order to bring together the discussion so far. Fix an
embedding $\Qbar \into \C$ and let $H$ be a multiset of integers of
cardinality $n>1$. Then the following sets should be in natural bijection.
One way to make precise the meaning of `natural' is to require that two
objects $M$ and $M'$ should correspond if the local L-factors $L_p(M,X)$ and
$L_p(M',X)$ are equal for all but finitely many $p$. Note that in each case
the factors $L_p(M,X)$ for all but finitely many $p$, completely determine
$M$.\vspace{3mm}

\begin{enumerate}
\item[(AF)]\label{af} Irreducible constituents $\pi$ of $\CA^\circ_H(GL_n(\Q) \backslash
GL_n(\A))$. \vspace{5mm}

\item[(LF)]\label{lf}
Near equivalence classes of irreducible, admissible $GL_n(\A^\infty)
\times (\gl_n,O(n))$-modules $\pi$ with the following properties. (We call
two $GL_n(\A^\infty) \times (\gl_n,O(n))$-modules, $\pi$ and $\pi'$ nearly
equivalent if $\pi_p \cong \pi_p'$ for all but finitely many $p$.)
\begin{enumerate}
\item $\pi_\infty$ has infinitesimal character $H$.
\item The central character $\psi_\pi$ is trivial on $\Q^{\times}
\subset \A^{\times}$.
\item For all characters $\psi:\A^\times/\Q^\times \R^\times_{>0}$ the
$L$-function $\Lambda(\pi \otimes (\psi \circ \det),s)$ converges in some right
half plane, has holomorphic continuation to the entire complex plane so that
it is bounded in vertical strips and satisfies the functional equation
\[  \Lambda(\pi \otimes \psi,s)=\epsilon(\pi \otimes \psi) N(\pi \otimes \psi
)^{-s} \Lambda(\pi^* \otimes \psi^{-1},1-s). \]
\item (See \cite{js} for an explanation of this condition.) There is a finite set of primes $S$ containing all primes $p$ for which
$\rec_p(\pi_p)$ is ramified, such that, writing $L(\pi_p,X) = \prod_{i=1}^n
(1-\alpha_{p,i}X)^{-1}$ for $p \not\in S$,
\[ \sum_{p\not\in S,i,j} \sum_{m=1}^{\infty} \alpha_{p,i}^m \alpha_{p,j}^{-m}
/mp^{ms} + \log(s-1) \]
is bounded as $s \ra 1$ from the right.
\end{enumerate}
In this case $L_p(\pi,X)=L(\pi_p,X)$. \vspace{5mm}

\begin{description}
\item[(lR)]\label{lr} (Fix $\iota: \Qbar_l \ra \C$.) Irreducible $l$-adic representations
\[ R:G_\Q \lra GL_n(\Qbar_l) \]
which are unramified at all but finitely many primes and for which
$R|_{G_{\Q_l}}$ is de Rham with $\HT(R|_{G_{\Q_l}})=H$. In this case
$L_p(R,X)=\iota L(\WD_p(R),X)$. \vspace{5mm}

\item[(WCS)]\label{cs} Irreducible weakly compatible systems of $l$-adic representations
$\CR$ with \linebreak $\HT(\CR)=H$. In this case $L_p(\CR,X)=L_p(\WD_p(\CR),X)$. \vspace{5mm}

\item[(GCS)]\label{gcs} Irreducible geometric strongly compatible systems of $l$-adic representations
$\CR$ with $\HT(\CR)=H$. In this case $L_p(\CR,X)=L_p(\WD_p(\CR),X)$.
\end{description}
\end{enumerate}

For $n=1$ we must drop the item $(LF)$, because it would need to be
modified to allow $L(\pi \otimes (\psi \circ \det),s)$ to have a simple
pole, while, in any case condition (LF) (b) would make the implication
$(LF) \Longrightarrow (AF)$ trivial. This being said,
in the case $n=1$ all the other four sets are known to be in natural bijection (see
\cite{serre:alr}). This basically follows because global class
field theory provides an isomorphism
\[ \Art: \A^{\times}/\Q^{\times} \R^\times_{>0} \liso G_\Q^\ab. \]

I would again like to stress how different are these various sorts
of objects and how surprising it is to me that there is any
relation between them. Items (AF) and (LF) both concern
representations of adele groups, but arising in rather different
settings: either from the theory of discrete subgroups of Lie
groups or from the theory of $L$-functions with functional
equation. Items (lR) and (WCS) arise from Galois theory and item
(GCS) arises from geometry.

So what do we know about the various relationships for $n>1$?

Not much. Trivially one has $(GCS) \Longrightarrow (WCS) \Longrightarrow (lR)$.
The passage $(AF) \Longrightarrow (LF)$ is OK by theorem \ref{gjt}. As discussed
in section \ref{s3} we have significant partial results in the directions
$(LF) \Longrightarrow (AF)$ and $(AF) \Longrightarrow (GCS)$, but both seem to
need new ideas. (Though I should stress that I am not really competent to discuss
converse theorems.)

One way to establish the equivalence of all five items would be to
complete the passages $(LF) \Longrightarrow (AF)$ and $(AF) \Longrightarrow
(GCS)$ and to establish the passage $(lR) \Longrightarrow (AF)$. It is these
implications which have received most study, though it should be pointed out
that in the function field case the equivalence of the analogous objects was
established by looking at the implications
\[ (lR) \Longrightarrow (LF) \Longrightarrow
(AF) \Longrightarrow (GCS). \]
(See \cite{laff}. It is the use of techniques from Grothendieck's $l$-adic
cohomology to prove the first of these implications which is most special to
function fields.) However it is striking that in the case of number
fields all known implications from items (lR), (WCS) or (GCS) to
(LF) go via (AF).

For the rest of this article we will concentrate on what still seems to be
the least understood problem: the passage from (lR) or (WCS) to (AF) or (LF).
Although the results we have are rather limited one should not underestimate
their power. Perhaps the most striking illustration of this is that the
lifting theorems discussed in section \ref{s4.2} (combined with earlier
work using base change and converse theorems) allowed Wiles \cite{wiles}
to finally prove Fermat's last theorem.

The discussion in the rest of this paper will of necessity be somewhat more
technical. In particular we will need to discuss automorphic forms, $l$-adic
representations and so on over general number fields (i.e. fields finite over $\Q$).
We will leave it to the reader's imagination exactly how
such a generalisation is made. In this connection we should remark that if
$L/K$ is a finite extension of number fields and if $R$ is a semi-simple
de Rham $l$-adic representation of $G_L$ which is unramified at all but
finitely many primes, then (see \cite{artin})
\[ L(R,s) = L(\Ind_{G_L}^{G_K} R,s) \]
(formally if the $L$-functions don't converge). In fact this is
true Euler factor by Euler factor and similar results hold for
conductors and $\epsilon$-factors (see \cite{tate:ntb}). This
observation can be extremely useful.

\subsection{Cyclic base change}\label{s4.1}

\vskip-5mm \hspace{5mm}

Suppose that $G$ is a group, $H$ a normal subgroup such that $G/H$ is cyclic
with generator of $\sigma$. It is an easy exercise that an irreducible
representation $r$ of $H$ extends to a representation of $G$ if and only if
$r^\sigma \cong r$ as representations of $H$. If one believes conjectures
\ref{cgfa} and \ref{cafg}, one might expect that if $L/K$ is a cyclic Galois
extension of number fields of prime order, if $\sigma$ generates
$\Gal(L/K)$ and if $\pi$ is a cuspidal automorphic representation of
$GL_n(\A_L)$ with $\pi \circ \sigma \cong \pi$, then there should be a
cuspidal automorphic representation $\Pi$ of $GL_n(\A_K)$, such that for all
places $w$ of $L$ we have $\rec_w(\pi_w)=\rec_{w|_K}(\Pi_{w|_K})|_{W_{L_w}}$.
This is indeed the case. For $n=1$ we have $\pi=\Pi \circ \norm_{L/K}$,
Langlands \cite{lang:bc} proved it for $n=2$ using the trace formula and
Arthur and Clozel \cite{ac} generalised his method to all $n$.

One drawback of this result is that if $v$ is a place of $K$ inert in
$L$ then there is no complete recipe for $\Pi_v$ in terms of $\pi$.
This can be surprisingly serious. It can however be alleviated, if we know
how to associate irreducible $l$-adic representations to both $\Pi$ and $\pi$.
Langlands used this to show that many two dimensional Artin representations
(i.e. $l$-adic representations with finite image) were automorphic (i.e.
associated to a cuspidal automorphic representation). In fact using additional
results from the theory of $L$-functions, particularly the converse theorem
for $GL_3$ (see section \ref{s4.4}), he and Tunnell (\cite{tunnell}) were
able to establish the automorphy of all continuous two dimensional Artin
representations with soluble image.

\subsection{Brauer's theorem}\label{s4.3}

\vskip-5mm \hspace{5mm}

The result I want to discuss is a result of Brauer \cite{brauer} about finite groups.
\begin{thm}[{\rm Brauer}] Suppose that $r$ is a representation of a finite group
$G$. Then there are soluble subgroups $H_i<G$, one dimensional representations
$\psi_i$ of $H_i$ and integers $n_i$ such that as virtual representations of
$G$ we have
\[ r = \sum_i n_i \Ind_{H_i}^G \psi_i. \]
\end{thm}

As Artin \cite{artin} had realised this theorem has the following immediate consequence.
(Indeed Brauer proved his theorem in response to Artin's work.)

\begin{cor} Let $\iota:\Qbar_l \ra \C$. Suppose that
\[ R:G_\Q \lra GL_n(\Qbar_l) \]
is an $l$-adic representation with {\em finite image}. Then the $L$-function
$L(\iota R, s)$ has meromorphic continuation to the entire complex plane and
satisfies the expected functional equation. \end{cor}

Artin's argument runs as follows. Let $G$ denote the image of $R$ and write
\[ R = \sum_i n_i \Ind_{H_i}^G \psi_i \]
as in Brauer's theorem. Let $L/\Q$ be the Galois extension with group $G$ cut
out by $R$ and let $K_i=L^{H_i}$. Then one has equalities
\[ \begin{array}{rcl} L(\iota R,s) &=& \prod_i L(\iota \Ind_{G_{K_i}}^{G_\Q}
\psi_i,s)^{n_i} \\ &=& \prod_i L( \iota \psi_i,s)^{n_i}. \end{array} \]
By class field theory for the fields $K_i$, the character $\psi_i$ is
automorphic on $GL_1(\A_{K_i})$ and so $L(\iota
\psi_i,s)$ has holomorphic continuation to the entire complex plane (except possibly for one simple pole if $\psi_i=1$) and
satisfies a functional equation. It follows that $L(\iota R,s)$ has
meromorphic continuation to the entire complex plane and satisfies a
functional equation.
The problem with this method is that some of the integers $n_i$ will usually
be negative so that one can only conclude the meromorphy of $L(\iota R,s)$,
not its holomorphy.

\subsection{Converse theorems}\label{s4.4}

\vskip-5mm \hspace{5mm}

As Cogdell and Piatetski-Shapiro point out, conjecture \ref{cct}
would have very important implications for Galois
representations. For instance the cases $n=2$ and $3$ played a key role
in the proof of the automorphy of two dimensional Artin representations
(see \ref{s4.1}). Conjecture \ref{cct} combined with Brauer's theorem
and a result of Jacquet and Shalika \cite{js} in fact implies that
many (all? - certainly those with soluble or perfect image) Artin
representations are automorphic. A similar argument shows that in many
other cases, in order to check the automorphy of an $l$-adic representation
of $G_\Q$, it suffices to do so after a finite base change. For instance one has
the following result.

{\em Assume conjecture \ref{cct}. Let $\iota:\Qbar_l \into \C$ and let $K/\Q$ be a
finite, totally real Galois extension. Suppose that $\Pi$ is a cuspidal
automorphic representation of $GL_n(\A_K)$ with infinitesimal character
corresponding to a multiset $H$ consisting of $n$ distinct integers. If $n>2$
also suppose that $\Pi_v$ is square integrable (i.e. $\rec_v(\Pi_v)$ is
indecomposable) for some finite place $v$ of $K$. Let
\[ R:G_\Q \lra GL_n(\Qbar_l) \]
be an $l$-adic representation such that $R \sim R^* \otimes \psi$ for some
character $\psi$ of $G_\Q$, and such that $R|_{G_{K}}$ is irreducible. Suppose
finally that $R|_{G_K}$ and $\Pi$ are associated, in the sense that, for all
but finitely many places $v$ of $K$, we have
\[ \iota L(\WD_v(R|_{G_K}),X) = L(\Pi_v,X). \]
Then there is a regular algebraic cuspidal automorphic representation $\pi$
of $GL_n(\A)$ associated to $R$ in the same sense.}

\subsection{Lifting theorems}\label{s4.2}

\vskip-5mm \hspace{5mm}

To describe this sort of theorem we first remark that if $R:G_\Q \ra
GL_n(\Qbar_l)$ is continuous then after conjugating $R$ by some element
of $GL_n(\Qbar_l)$ we may assume that the image of $R$ is contained in
$GL_n(\CO_{\Qbar_l})$ and so reducing we obtain a continuous representation
\[ \barR: G_\Q \lra GL_n(\Fbar_l). \]

The lifting theorems I have in mind are results of the general form if
$R$ and $R'$ are $l$-adic representations of $G_\Q$ with $R'$ automorphic
and if $\barR=\barR'$ then $R$ is also automorphic. Very roughly speaking
the technique (pioneered by Wiles \cite{wiles} and completed by the author and
Wiles \cite{tw}) is to show that $R \bmod l^r$ arises from automorphic forms
for all $r$ by induction on $r$. As $\ker(GL_n(\Z/l^r\Z) \onto GL_n(\Z/l^{r-1}
\Z))$ is an abelian group one is led to questions of class field theory
and Galois cohomology.

I should stress that such theorems are presently available only in very limited
situations. I do not have the space to describe the exact limitations, which
are rather technical, but
the sort of restrictions that are common are as follows.
\begin{enumerate}
\item If $R:G_\Q \ra GL(V)$ then there should be a character $\mu:G_\Q \ra
GL_n(\Qbar_l)$ and a non-degenerate bilinear form $(\,\,\, ,\,\,\, )$ on $V$
such that
\begin{itemize}
\item $(R(\sigma)v_1,R(\sigma)v_2) = \mu(\sigma)(v_1,v_2)$ and
\item $(v_2,v_1) = \mu(c) (v_1,v_2)$.
\end{itemize}
(This seems to be essential for the method of \cite{tw}.)

\item $R$ should be de Rham with distinct Hodge-Tate numbers. (This again
seems essential to the method of \cite{tw}, but see \cite{bt}.)

\item Either $R$ and $R'$ should be ordinary (i.e. their restrictions to
$G_{\Q_l}$ should be contained in a Borel subgroup); or $R$ and $R'$ should be crytsalline (not just de Rham) at $l$ with
the same Hodge-Tate numbers and $l$ should be large compared with the
differences of elements of $\HT(R)$. (The problems here are connected with
the need for an integral Fontaine theory, but they are not simply technical
problems. There are some complicated results pushing back this restriction in
isolated cases, see \cite{cdt}, \cite{bcdt}, \cite{savitt}, but so far our
understanding is very limited.)

\item The image of $\barR$ should not be too small (e.g. should be irreducible
when restricted to $\Q(e^{2 \pi i/l})$),
though in the case $n=2$ there is beautiful work of Skinner and Wiles
(\cite{sw1} and \cite{sw2}) dispensing with this criterion, which this author
has unfortunately not fully understood.

\end{enumerate}

In addition, all the published work is for the case $n=2$. However
there is ongoing work of a number of people attempting to dispense
with this assumption. Using a very important insight of Diamond
\cite{diam}, the author, together with L.Clozel and M.Harris, has
generalised to all $n$ the so called minimal case (originally
treated in \cite{tw}) where $R$ is no more ramified than $\barR$.
One would hope to be able to deduce the non-minimal case from
this, as Wiles did in \cite{wiles} for $n=2$. In this regard one
should note the work of Skinner and Wiles \cite{swe} and the work
of Mann \cite{russ}. However there seems to be one missing
ingredient, the analogue of the ubiquitous Ihara lemma, see lemma
3.2 of \cite{ihara} (and also theorem 4.1 of \cite{rib2}). As this
seems to be an important question, but one which lies in the
theory of discrete subgroups of Lie groups, let us take the
trouble to formulate it, in the hope that an expert may be able to
prove it. It should be remarked that there are a number of
possible formulations, which are not completely equivalent and any
of which would seem to suffice. We choose to present one which has
the virtue of being relatively simple to state.

\begin{conjecture} Suppose
that $G/\Q$ is a unitary group which becomes an inner form of $GL_n$ over
an imaginary quadratic field $E$. Suppose that $G(\R)$ is compact. Let
$l$ be a prime which one may assume is large compared to $n$. Let
$p_1$ and $p_2$ be distinct primes different from $l$ with $G(\Q_{p_1}) \cong GL_n(\Q_{p_1})$ and
$G(\Q_{p_2}) \cong GL_n(\Q_{p_2})$. Let $U$ be an open compact subgroup of
$G(\A^{p_1,p_2})$ and consider the representation of $GL_n(\Q_{p_1}) \times
GL_n(\Q_{p_2})$ on the space $C^\infty(G(\Q)\backslash G(\A) / U, \Fbar_l)$
of locally constant $\Fbar_l$-valued functions on
\[ G(\Q) \backslash G(\A) /U = (G(\Q) \cap U) \backslash (GL_n(\Q_{p_1}) \times
GL_n(\Q_{p_2}). \]
(Note that $G(\Q) \cap U$ is a discrete cocompact subgroup of $GL_n(\Q_{p_1})
\times GL_n(\Q_{p_2})$.) Suppose that $\pi_1 \otimes \pi_2$ is an
irreducible {\em sub-representation} of $C^\infty(G(\Q)\backslash G(\A) / U,
\Fbar_l)$ with $\pi_1$ generic. Then $\pi_2$ is also generic.
\end{conjecture}

The most serious problem with applying such lifting theorems to
prove an $l$-adic representation $R$ is automorphic is the need to
find some way to show that $\barR$ is automorphic. The main
success of lifting theorems to date, has been to show that if $E$
is an elliptic curve over the rationals then $H^1(E(\E),\Qbar_l)$
is automorphic, so that $E$ is a factor of the Jacobian of a
modular curve and the $L$-function $L(E,s)$ is an entire function
satisfying the expected functional equation (\cite{wiles},
\cite{tw},\cite{bcdt}). This was possible because $GL_2(\Z_3)$
happens to be a pro-soluble group and there is a homomorphism
$GL_2(\F_3) \lra GL_2(\Z_3)$ splitting the reduction map. The
Artin representation
\[ G_\Q \lra GL(H^1(E(\C),\F_3)) \lra GL_2(\Z_3) \]
is automorphic by the Langlands-Tunnell theorem alluded to in section
\ref{s4.1}.

\subsection{Other techniques?}

\vskip-5mm \hspace{5mm}

I would like to discuss one other technique which has been some
help if $n=2$ and may be helpful more generally. We will restrict
our attention here to the case $n=2$ and $\det R(c)=-1$. We have
said that the principal problem with lifting theorems for proving
an $l$-adic representation $R:G_{\Q} \ra GL_2(\Qbar_l)$ is
automorphic is that one needs to know that $\barR$ is automorphic.
This seems to be a very hard problem. Nonetheless one can often
show that $\barR$ becomes automorphic over some Galois totally
real field $K/\Q$. (Because $K$ is totally real, if $\barR (G_\Q)
\supset SL_2(\F_l)$ and $l>3$ then $\barR (G_K) \supset
SL_2(\F_l)$. So this `potential automorphy' is far from vacuous).
The way one does this is to look for an abelian variety $A/K$ with
multiplication by a number field $F$ with $[F:\Q]=\dim A$, and
such that $\barR$ is realised on $H^1(A(\C),\F_l)[\lambda]$ for
some prime $\lambda|l$, while for some prime $\lambda'|l' \neq l$
the image of $G_K$ on $H^1(A(\C),\F_{l'})[\lambda']$ is soluble.
One then argues that $H^1(A(\C),\F_{l'})[\lambda']$ is
automorphic, hence by a lifting theorem $H^1(A(\C),\Q_{l'})
\otimes_{F_{l'}} F_{\lambda'}$ is automorphic, so that
(tautologically) $H^1(A(\C),\F_l)[\lambda]$ is also automorphic,
and hence, by another lifting theorem, $R|_{G_K}$ is automorphic.
One needs $K$ to be totally real, as over general number fields
there seems to be no hope of proving lifting theorems, or even of
attaching $l$-adic representations to automorphic forms. In
practice, because of various limitations in the lifting theorems
one uses, one needs to impose some conditions on the behaviour of
a few primes, like $l$, in $K$ and some other conditions on $A$.
The problem of finding a suitable $A$ over a totally real field
$K$, comes down to finding a $K$-point on a twisted Hilbert
modular variety. This is possible because we are free to choose
$K$, the only restriction being that $K$ is totally real and
certain small primes (almost) split completely in $K$. To do this,
one has the following relatively easy result.

\begin{prop}[{\rm \cite{bailly},\cite{pop}}] Suppose that $X/\Q$ is a smooth geometrically
irreducible variety. Let $S$ be a finite set of places of $\Q$ and suppose that
$X$ has a point over the completion of $\Q$ at each place in $S$. Let
$\Q_S$ be the maximal extension of $\Q$ in which all places in $S$
split completely (e.g. $\Q_{\{ \infty \}}$ is the maximal totally
real field). Then $X$ has a $\Q_S$-point. \end{prop}

In this regard it would have extremely important consequences if, in the
previous proposition, one could replace $\Q_S$ by $\Q_S^\sol$, the maximal
soluble extension of $\Q$ in which all places in $S$ split completely.
I do not know if it is reasonable to expect this.

Using this method one can, for instance, prove the following result.

\begin{thm}[{\rm \cite{mero}}] \label{mer}
Suppose that $\CR$ is an irreducible weakly compatible system of two
dimensional $l$-adic representations with $\HT(\CR)=\{ n_1,n_2 \}$ where
$n_1 \neq n_2$. Suppose also that $\det R_{l,\iota}(c)=-1$ for one
(and hence for all) pairs $(l,\iota)$. Then there is a Galois totally
real field $K/\Q$ and a cuspidal automorphic representation $\pi$ of
$GL_2(\A_K)$ such that
\begin{itemize}
    \item for all $v|\infty$, $\pi_v$ has infinitesimal character
    $H$, and
    \item for all $(l,\iota)$ and for all finite places $v \ndiv l$
    of $K$ we have
    \[ \rec_v(\pi_v) = \WD_v(R_{l,\iota}|_{G_K})^\semis. \]
\end{itemize}
In particular $\CR$ is pure of weight $(n_1+n_2)/2$. Moreover $\CR$ is
strongly compatible and $L(\iota \CR,s)$ has meromorphic
continuation to the entire complex plane and satisfies the expected
functional equation.
\end{thm}

The last sentence of this theorem results from the first part and Brauer's
theorem. We remark that conjecture \ref{cct} would imply that this
theorem could be improved to assert the automorphy of $\CR$ over $\Q$.

\label{lastpage}

\end{document}